\input amstex
\documentstyle{amsppt}
\magnification=\magstep1
\define \my{\bf}
\define\po{\parindent 0pt}
\define \mytitle{\po \m \my}
\define \p{{\po \it \m Proof. }}
\define \m{\medskip}
\nologo
\TagsOnRight
\define \ph#1{\phantom {#1}}
\define \T#1{\widetilde {#1}}

\define \IM{\operatorname {Im}}
\define \cat{\operatorname{cat}}
\define \swgt{\operatorname{swgt}}
\define \crit{\operatorname{crit}}
\define \Crit{\operatorname{Crit}}
\define \Ext{\operatorname{Ext}}
\define \Hom{\operatorname{Hom}}
\define \Rest{\operatorname{Rest}}
\define \sgrad{\operatorname{sgrad}}
\define \Fix{\operatorname{Fix}}

\define \eps{\varepsilon}
\def \sirc{{\raise0.2ex \hbox{$\scriptstyle \circ$}}}
\define \mul{{\raise0.2ex \hbox{$\scriptstyle \bullet$}}}
\topmatter
\title{ON STRICT CATEGORY WEIGHT AND  THE ARNOLD CONJECTURE}\endtitle
\author{Yuli B. Rudyak}\endauthor
\leftheadtext{Yuli B. Rudyak}
\rightheadtext{Arnold conjecture}
\address{Universit\"at--GH Siegen, FB6/Mathematik, 57068 Siegen,
Germany}
\endaddress 
\email{rudyak\@mathematik.uni-siegen.de}\endemail
\abstract In [R2] and [RO] the Arnold conjecture for symplectic  
manifolds $(M,\omega)$ with $\pi_2(M)=0$ was proved. This proof used  
surgery and cobordism theory. Here we give a purely cohomological  
proof of this result.
\endabstract
\endtopmatter
\head{Introduction}\endhead
Given a smooth (=$C^{\infty}$) manifold $M$, we set
$
\Crit M:= \min\{\crit f\}
$
 where $f$ runs over all smooth functions $M \to \Bbb R$. 
\m Let $(M,\omega)$ be a symplectic manifold. Given a smooth  
function $f: M \to \Bbb R$, let $\sgrad f$ denote the symplectic  
gradient of $f$, i.e., the vector field defined as follows:
$$
\omega(\sgrad f, \xi)=-df(\xi)
$$
for every vector field $\xi$.
\m In [A] Arnold proposed the following remarkable conjecture. Let  
$(M,\omega)$ be a closed symplectic manifold, and let $H: M \times   
\Bbb R \to \Bbb R$ be a smooth function such that $H(m,t)=H(m,t+1)$  
for every $m \in M, t\in \Bbb R$. We define $H_t: M \to \Bbb R$ by  
setting $H_t(m):= H(m,t)$. Consider the time-dependent differential  
equation
$$
\dot x= \sgrad H_t(x(t)).
\tag{*}
$$
The Arnold conjecture claims that the number of 1-periodic  
solutions of (*) is at least $\Crit M$. 
\m This conjecture admits another interpretation. The equation (*)  
yields a family $\varphi_t: M \to M, t\in \Bbb R$, where, for every  
$x\in M$,  $\varphi_t(x)$ is the integral curve of (*). A {\it  
Hamiltonian symplectomorphism} is a diffeomorphism $\phi: M \to M$  
which has the form $\phi=\varphi_1$ for some function $H: M \times  
\Bbb R \to \Bbb R$ as above. So, the Arnold conjecture can be  
refolmulated as follows:
$$
\Fix \phi \geq \Crit M
$$
for every Hamiltonian symplectomorphism $\phi$, where $\Fix \phi$  
denotes the number of fixed points of $\phi$. 
\m This conjecture was proved for many special cases, see [MS] and  
[HZ] for a survey. Here we notice the following result of Floer~[Fl]  
and Hofer~[H]: the number $\Fix \phi$ can be estimated from below  
by the cup-lenght of $M$. So, here we have a weak form of the Arnold  
conjecture.
\m In [R2] and [RO] the Arnold conjecture was proved for every  
closed connected symplectic manifold $(M,\omega)$ with  
$\omega|_{\pi_2(M)}=0=c_1|_{\pi_2(M)}$. In greater detail, in [R2]  
the conjecture was proved under the additional condition $\cat  
M=\dim M$, and it was proved in [RO] that the condition  
$\omega|_{\pi_2(M)}=0$ implies the condition $\cat M=\dim M$.  
Because of the last result, it turns out to be that $\Crit M=2n+1$  
provided $\omega|_{\pi_2(M)}=0$, and actually we have the inequality  
$\Fix \phi \geq 2n+1$.
\par The proof of the Arnold conjecture in [R] uses surgery and  
cobordism theory. Here we give another proof of the Arnold  
conjecture (under the same restriction  
$\omega|_{\pi_2(M)}=0=c_1|_{\pi_2(M)}$). This proof uses the  
ordinary cohomology $H^*$ only; probably, it is more convenient for  
people which work in the area of dynamical systems and are not very  
familiar with extraordinary cohomology. The main line of the proof  
follows Rudyak--Oprea~[RO], but here we use the strict category  
weight instead the category weight.
{\mytitle Remarks}. 1. Hofer and Zehnder~[HZ, p.250] mentioned that  
the Theorem 2.2 below is true without the restriction  
$c_1|_{\pi_2(M)}=0$. Because of this, the Arnold conjecture turns  
out to be valid for all closed connected symplectic manifolds with  
$\omega|_{\pi_2(M)}=0$.  
\par 2. Actually, in Theorem 2.2 below the number $\Fix \phi$ is  
the number of 1-periodic solutions of the equation $(*)$, while  
$\Rest \Phi$ is the number of contractble 1-periodic solutions of  
this equation. So, here (as well as in [R2] and [RO]) it is proved  
that the number of contractible 1-periodic solutions of $(*)$ is at  
least $2n+1$.
\m The paper is organized as follows. In \S 1 we discuss strict  
category weight, in \S 2 we use Floer's results in order to reduce  
the Arnold conjecture to a certain topological problem, in \S 3 we  
prove main results, in Appendix we discuss an analog of the Arnold  
conjecture for locally Hamiltonian symplectomorphisms.
\m The cohomology group $H^n(X;G)$ is always defined to be the  
Alexander--Spanier cohomology group with coefficient group $G$, see  
[M] or [S] for the details.
\par We reserve the term ``map'' for continuous functions.
\par ``Connected'' always means path connected.
\head \S 1. Strict category weight \endhead
{\mytitle 1.1. Definition {\rm ([LS], [Fox], [F], [BG])}.} Given a  
map $\varphi: A \to X$, we define the {\it Lusternik--Schnirelmann  
category} $\cat \varphi $ of $\varphi$ to be the minimal number $k$  
with the following property: $A$ can be covered by open sets  
$A_1,\ldots, A_{k+1}$ such that $\varphi|A_i$ is null-homotopic for  
every $i$. Furthrmore, we define the {\it Lusternik--Schnirelmann  
category} $\cat X$ of a space $X$ by setting  $\cat X:= \cat 1_X$.
\proclaim{1.2. Proposition {\rm ([BG])}} {\rm (i)}  For every  
diagram $A @>\varphi>> Y @>f>> X$ we  have $\cat f\varphi \leq  
\min\{\cat \varphi, \cat f\}$. In  particular, $\cat f\leq \min  
\{\cat X, \cat Y\}$. \par 
{\rm (ii)} If $\varphi \simeq \psi: A \to X$ then $\cat \varphi  
=\cat \psi$. \par 
{\rm (iii)} If $h: Y \to X$ is a homotopy equivalence then $\cat  
\varphi=\cat h \varphi$ for every $\varphi: A \to X$.
\qed
\endproclaim
\m Given a connected pointed space $X$, let $\eps: S\Omega X \to X$  
be the map adjoint to $1_{\Omega X}$, here $\Omega X$ is the loop  
space of $X$ and $S$ denotes the suspension, see e.g. [Sw].
\proclaim{1.3. Theorem {\rm ([Sv, Theorems 3, 19$'$ and 21])}} Let  
$\varphi: A \to X$ be a map of connected Hausdorff paracompact  
spaces. Then $\cat \varphi<2$ iff there is a map $\psi: A \to  
S\Omega X$ such that $\eps\psi=\varphi$. 
\qed
\endproclaim
{\mytitle 1.4. Definition {\rm ([R1])}.} Let $X$ be a Hausdorff  
paracompact space, and let $u \in H^q(X;G)$ be an arbitrary element.  
We define the {\it strict category weight} of $u$ (denoted by  
$\swgt u$) by setting
$$
\swgt u=\sup\{k\bigm|\varphi^*u=0 \text{ for every map $\varphi: A  
\to X$ with } \cat \varphi<k\}
$$
where $A$ runs over all Hausdorff paracompact spaces.
\m We use the term ``strict category weight'', since the term  
``category weight'' is already used (introduced) by  
Fadell--Hussein~[FH]. Concerning the relation between category  
weight and strict category weight, see~[R1].
\m We remark that $\swgt u={\infty}$ if $u=0$. \proclaim{1.5.  
Theorem} Let $X$ and $Y$ be two Hausdorff paracompact spaces. Then  
for every $u\in H^*(X)$ the following hold: \par 
{\rm (i)} for every map $f: Y \to X$ we have $\cat f \geq \swgt u$  
provided $f^*u \neq 0$. Furthermore, if $X$ is connected then $\swgt  
u \geq 1$ whenever $u\in \T H^*(X)$;\par 
{\rm (ii)} for every map $f: Y \to X$ we have $\swgt f^*u \geq  
\swgt u$; \par
{\rm (iii)} for every $u,v\in H^*(X)$ we have $\swgt(uv)\geq \swgt  
u + \swgt v$. 
\endproclaim
\p (i) This follows from the definition of $\swgt$.
\par (ii) This follows from 1.2(i).
\par (iii) Let $\swgt u=k,\ \swgt v=l$ with $k,l<\infty$. Given $f:  
A \to X$ with $\cat f < k+l$, we prove that $f^*(uv)=0$. Indeed,  
$\cat f<k+l$, and so $A=A_1\cup \cdots\cup A_{k+l}$ where each $A_i$  
is open in $A$ and $f|A_i$ is null-homotopic. Set $B=A_1\cup  
\cdots\cup A_k$ and $C=A_{k+1}\cup \cdots \cup A_{k+l}$. Then $\cat  
f|B<k$ and $\cat f|C<l$. Hence $f^*u|B=0=f^*v|C$, and so  
$f^*(uv)|(B\cup C)=0$. i.e., $f^*(uv)=0$.\par
The case of infinite category weight is leaved to the reader.
\qed
\head \S 2. Floer's reduction and related things.\endhead
{\mytitle 2.1. Recollection.} A {\it flow} on a topological space  
$X$ is a family $\Phi=\{\varphi_t\}, t\in \Bbb R$ where each  
$\varphi_t: X \to X$ is a self-homeomorphism and  
$\varphi_s\varphi_t=\varphi_{s+t}$ for every $s,t\in \Bbb R$ (notice  
that this implies $\varphi_0=1_X$). 
\par A flow is called {\it continuous} if the function $X \times  
\Bbb R \to X, (x,t)\mapsto \varphi_t(x)$ is continuous.
\par A point $x\in X$ is called a {\it rest point} of $\Phi$ if  
$\varphi_t(x)=x$ for every $t\in \Bbb R$. We denote by $\Rest \Phi$  
the number of rest points of $\Phi$.
\par A continuous flow $\Phi=\{\varphi_t\}$ is called {\it  
gradient-like} if there exists a continuous (Lyapunov) function $F:  
X \to \Bbb R$ with the following property: for every $x\in X$ we  
have $F(\varphi_t(x))< F(\varphi_s(x))$ whenever $t>s$ and $x$ is  
not a rest point of $\Phi$. 
\m The following Theorem can be found in [Fl, Th. 7], cf. also [HZ].
\proclaim{2.2. Theorem} Let $(M,\omega)$ be a closed connected  
symplectic manifold such that  
$\omega|_{\pi_2(M)}=0=c_1|_{\pi_2(M)}$, and let $\phi: M \to M$ be a  
Hamiltonian symplectomorphism. Then there exists a map $f: X \to M$  
with the following properties: \par 
{\rm (i)} $X$ is a compact metric space;\par 
{\rm (ii)} $X$ possesses a continuous gradient-like flow $\Phi$  
such that $\Rest \Phi \leq \Fix \phi$; \par 
{\rm (iii)} The homomorphism $f^*: H^n(M;G) \to H^n(X;G)$ is a  
monomorphism for every coefficient group $G$.
\qed
\endproclaim
The following theorem (of Lusternik--Schnirelmann type) is proved  
in [R2]. 
\proclaim{2.3. Theorem} Let $\Phi$ be a continuous gradient-like  
flow on a compact metric space $X$, let $Y$ be a Hausdorff space  
which can be covered by open and contractible in $Y$ subspaces, and  
let $f: X \to Y$ be a map. Then
$$
\Rest \Phi \geq 1+\cat f. \qed
$$
\endproclaim
We need also the following well-known fact which follows from [LS]  
and [T].
\proclaim{2.4. Theorem} For every closed smooth manifold $M$ we have 
$$
1+\cat M \leq \Crit M \leq 1+\dim M. \qed
$$
\endproclaim
\head \S 3. Proof of the Arnold conjecture\endhead
\proclaim{3.1. Theorem {\rm (cf. [FH], [RO], [St])}} Let $\pi$ be a  
discrete group. Then for every $u\in H^k(K(\pi,1);G)$ with $k>1$ we  
have $\swgt u \geq 2$.
\endproclaim
(Actually, Strom~[St] proved that $\swgt u\geq k$. Moreover, it is  
easy to see that $\swgt u \leq k$ provided $u \neq 0$, and so $\swgt  
u=k$ if $u \neq 0$.)
\p Because of 1.3, it suffices to prove that $\eps^*u=0$ where  
$\eps$ is a map from 1.3 and $\eps^*; H^*(K(\pi,1);G) \to  
H^*(S\Omega K(\pi,1);G)$ is the induced homomorphism. But $\Omega  
K(\pi,1)$ is homotopy equivalent to the discrete space $\pi$, and so  
$S\Omega K(\pi,1)$ is homotopy equivalent to a wedge of circles.  
Hence, $H^i(K(\pi,1(;G)=0$ for $k>1$, and thus $\eps^*u=0$.
\qed 
\proclaim{3.2. Theorem {\rm (cf. [RO])}} Let $Y$ be a connected  
finite $CW$-space, and let $y\in H^2(Y;G)$ be such that  
$y|\pi_2(Y)=0$. Then $\swgt y\geq2$.
\endproclaim
\p Let $\pi=\pi_1(Y)$, and let $g: Y \to K(\pi,1)$ be a map which  
induces an isomorphism of fundamental groups. First, we prove that 
$$
y\in \IM \{g^*: H^2(K(\pi,1);R) \to H^2(Y;R)\}.
$$
 Indeed, since $Y$ is a finite $CW$-space, its singular cohomology  
coincides with $H^*$, and so we have the universal coefficient  
sequence
$$
0 \to \Ext (H_1(Y), R) @>>> H^2(Y;R) @>l>> \Hom(H_2(Y),R) \to 0.
$$
On the other hand, there is a Hopf exact sequence
$$
\pi_2(Y) \to H_2(Y) \to H_2(K(\pi,1))  \to 0,
$$
and so we have the following commutative diagram with exact rows  
and colomn:
$$
\CD
\Ext (H_1(K), R) @>>> H^2(K);R) @>>>  \Hom(H_2(K),R))@>>>  0\\
@Vg'V\cong V @VVg^*V @VVg''V @.\\
 \Ext (H_1(Y), R) @>>>  H^2(Y;R) @>l>> \Hom(H_2(Y),R) @>>> 0\\
@. @. @VVV  @.\\
@. @.  \Hom(\pi_2(Y),R) @.
\endCD
$$
where $K$ denotes $K(\pi,1)$. Now, since $y|\pi_2(Y)=0$, we  
conclude that $l(y)\in \IM g''$. Since $g'$ is an isomorphism, an  
easy diagram hunting shows that $y\in \IM g^*$. 
\par Thus, by 3.1 and 1.5(ii), $\swgt y\geq2$.
\qed
\proclaim{3.3. Corollary} Let $Y$ be a connected finite $CW$-space,  
let $R$ be a commutative ring, let $y\in H^2(Y;R)$ be such that  
$y|\pi_2(Y)=0$,  and let $X$ be a Hausdorff paracompact space. 
If $f: X \to Y$ is a map with $f^*(y^n) \neq 0$, then $\cat f \geq 2n$.
\endproclaim
\p By 1.5,
$$
\cat f\geq \swgt y^n \geq n\swgt y \geq 2n.
\qed 
$$
\proclaim{3.4. Corollary {\rm ([RO])}} Let $(M^{2n},\omega)$ be a  
closed connected symplectic manifold with $\omega|_{\pi_2(M)}=0$.  
Then $\cat M=2n$ and $\Crit M=2n+1$.
\endproclaim
\p Since $\omega^n\neq 0$, we conclude that, by 3.3, $\cat M=\cat  
1_M \geq 2n$.
So, $\cat M=2n$ since $\cat M \leq \dim M$. The second equality  
follows from 2.4.
\qed
\proclaim{3.5. Corollary} Let $(M^{2n},\omega)$ be a closed  
symplectic manifold with $\omega|_{\pi_2(M)}=0=c_1|_{\pi_2(M)}$, and  
let $\phi: M \to M$ be a Hamiltonian symplectomorphism. Then $\Fix  
\phi \geq 2n+1$. 
In  particular, the Arnold conjecture holds for $(M,\omega)$.
\endproclaim
\p Let $f: X \to M$ and $\Phi$ be a map as in 2.2. Since every  
closed connected smooth manifold is a finite polyhedron, and since  
$\omega^n$ yields a non-trivial cohomology class in $H^*(X;\Bbb R)$,  
we conclude that, by 3.3, $\cat f \geq 2n$. Now, by 2.2 and 2.3,
$$
\Fix \phi \geq \Rest \Phi \geq 1+\cat f \geq 2n+1.
$$
Thus, because of 3.4, the Arnold conjecture holds for $(M,\omega)$.
\qed
\head{References}\endhead
\hyphenation{To-pol-ogy}
\halign{{\bf #\ }\hfil & \vtop{\parindent0pt
\hsize=31.1em
\hangindent0em\strut#\strut}\cr
[A] & V. I. Arnold:\it Mathematical Methods in Classical Mechanics.  
\rm Springer, Berlin Heidelberg New York 1989\cr
[BG] & I. Berstein, T. Ganea:  \it The category of a map and of a  
cohomology class.  \rm  Fund. Math. {\bf 50} (1961/2), 265--279\cr
[FH] & E. Fadell, S. Husseini: Category weight and Steenrod  
operations, Boletin de la Sociedad Matem\'atica Mexicana, {\bf 37}  
(1992), 151--161\cr
[F] & A. I. Fet: \it A connection between the topological  
properties and the number of extremals on a manifold  \rm (Russian).  
 Doklady AN SSSR, {\bf 88} (1953), 415--417 \cr 
[Fl] & A. Floer: \it Symplectic fixed points and holomorphic  
spheres, \rm Commun. Math. Phys. {\bf 120} (1989), 575--611\cr
[Fox] & R. Fox: \it On the Lusternik--Schnirelmann category\rm ,  
Ann. of Math. {\bf 42} (1941), 333-370\cr 
[H] & H. Hofer: \it Lusternik--Schnirelmann theory for Lagrangian  
intersections. \rm Annales de l'inst. Henri Poincar\'e -- analyse  
nonlineare, {\bf  5} (1988), 465--499\cr
[HZ] & H. Hofer, E. Zehnder: {\it Symplectic Invariants and  
Hamiltonian Dynamics}, Birkh\"auser, Basel, 1994\cr
[LS] & L. A. Lusternik, L. G. Schnirelmann,   {\it Methodes  
topologiques dans le probl\`emes variationels}. Hermann, Paris 1934  
\cr
[M] & W. Massey: {\it Homology and cohomology theory}, Marcel  
Dekker, INC, New York and Basel, 1978\cr
[MS] & D. McDuff, D. Salamon: \it Introduction to Symplectic  
Topology, \rm Clarendon Press, Oxford 1995\cr
[R1] & Yu. B. Rudyak: \it On category weight and its applications,  
\rm to appear in Topology, (1998)\cr
[R2] & Yu. B. Rudyak: \it On analytical applications of stable  
homotopy $($the Arnold conjecture, critical points$)$ \rm to appear  
in Math. Z., available as dg-ga/9708008\cr
[RO] & Yu. B. Rudyak, J. Oprea: \it On the Lusternik--Schnirelmann  
Category of Symplectic  
Manifolds and the Arnold Conjecture, \rm to appear in Math. Z.,  
available as dg-ga/9708007 \cr
[S] & E. H. Spanier: \it Algebraic Topology. \rm McGraw-Hill, New  
York 1966.\cr
[St] & J. Strom: {\it Essential category weight and classifying  
spaces}, Preprint Univ. of Madison-Wisconsin, 1997.\cr
[Sv] & A. Svarc: \it The genus of a fiber space. \rm Amer. Math.  
Soc. Translations {\bf 55} (1966), 49--140\cr
[Sw] & R. Switzer: \it Algebraic Topology -- Homotopy and Homology.  
\rm Springer, Berlin Heidelberg New York (1975)\cr
[T] & F. Takens: \it The minimal number of crirical points of a  
function on a compact manifold and the Lusternik--Schnirelmann  
category, \rm Invent. Math. {\bf 6} (1968), 197--244 \cr }
\enddocument